\documentclass{ws-rv975x65}
\usepackage[square]{ws-rv-van}             % numbered citation/references
\usepackage{tikz}

\makeindex
\mathchardef\ordinarycolon\mathcode`\:
  \mathcode`\:=\string"8000
  \begingroup \catcode`\:=\active
    \gdef:{\mathrel{\mathop\ordinarycolon}}
  \endgroup
\newtheorem{construction}[theorem]{Construction}

\begin{document}

\chapter{Coding theory and algebraic combinatorics\label{Huber_CodingAlgComb}}\index{coding theory}\index{algebraic combinatorics}\index{combinatorics}

\author[Michael Huber]{Michael Huber\footnote{The author gratefully acknowledges support by the Deutsche
Forschungsgemeinschaft (DFG).}}

\address{Institut f\"{u}r Mathematik, Technische Universit\"{a}t Berlin,\\
Stra{\ss}e des 17.~Juni~136, D-10623 Berlin, Germany, \\
mhuber@math.tu-berlin.de}

\begin{abstract}
This chapter introduces and elaborates on the fruitful interplay of
coding theory and algebraic combinatorics, with most of the focus on
the interaction of codes with combinatorial designs, finite
geometries, simple groups, sphere packings, kissing numbers,
lattices, and association schemes. In particular, special interest
is devoted to the relationship between codes and combinatorial
designs. We describe and recapitulate important results in the
development of the state of the art. In addition, we give
illustrative examples and constructions, and highlight recent
advances. Finally, we provide a collection of significant open
problems and challenges concerning future research.
\end{abstract}

\body

%%% ----------------------------------------------------------------------

\section{Introduction}\label{intro}

The classical publications ``A mathematical theory of
communication'' by C.~E.~Shannon~\cite{shan48} and ``Error detecting
and error correcting codes'' by R.~W.~Hamming~\cite{ham50} gave
birth to the twin disciplines of information theory and coding
theory. Since their inceptions the interactions of information and
coding theory with many mathematical branches have continually
deepened. This is in particular true for the close connection
between coding theory and algebraic combinatorics.

%\smallskip

This chapter introduces and elaborates on this fruitful interplay of
coding theory and algebraic combinatorics, with most of the focus on
the interaction of codes with combinatorial designs, finite
geometries, simple groups, sphere packings, kissing numbers,
lattices, and association schemes. In particular, special interest
is devoted to the relationship between codes and combinatorial
designs. Since we do not assume the reader is familiar with the
theory of combinatorial designs, an accessible and reasonably
self-contained exposition is provided. Subsequently, we describe and
recapitulate important results in the development of the state of
the art, provide illustrative examples and constructions, and
highlight recent advances. Furthermore, we give a collection of
significant open problems and challenges concerning future research.

\smallskip

The chapter is organized as follows. In Sec.~\ref{back}, we give a
brief account of basic notions of algebraic coding theory.
Section~\ref{pract} consists of the main part of the chapter: After
an introduction to finite projective planes and combinatorial
designs, a subsection on basic connections between codes and
combinatorial designs follows. The next subsection is on perfect
codes and designs, and addresses further related concepts.
Subsection~\ref{assmatt_thm} deals with the classical Assmus-Mattson
Theorem and various analogues. A subsection on codes and finite
geometries follows the discussion on the non-existence of a
projective plane of order $10$. In
Subsection~\ref{golay_witt_mathieu}, interrelations between the
Golay codes, the Mathieu-Witt designs, and the Mathieu groups are
studied. Subsection~\ref{golay_lattice} deals with the Golay codes
and the Leech lattice, as well as recent milestones concerning
kissing numbers and sphere packings. The last topic of this section
considers codes and association schemes. The chapter concludes with
sections on directions for further research as well as conclusions
and exercises.

%%% ----------------------------------------------------------------------

\section{Background}\label{back}

For our further purposes, we give a short account of basic notions
of algebraic coding theory. For additional information on the
subject of algebraic coding theory, the reader is referred
to~\cite{berl68,cali91,hill86,hupl_handb98,hupl03,kerb06,lint95,lint99,pewa72,roth06,wisl77}.
For some historical notes on its origins, see~\cite{thomp83}
and~\cite[Chap.\,1]{hupl_handb98}, as well as~\cite{cald98} for a
historical survey on coding theory and information theory.

\smallskip

We denote by $\mathbb{F}^n$ the set of all $n$-tuples from a
$q$-symbol alphabet. If $q$ is a prime power, we take the finite
field $\mathbb{F}=\mathbb{F}_q$ with $q$ elements, and interpret
$\mathbb{F}^n$ as an \mbox{$n$-dimensional} vector space
$\mathbb{F}_q^n$ over $\mathbb{F}_q$. The elements of $\mathbb{F}^n$
are called \emph{vectors} (or \emph{words})\index{word} and will be
denoted by bold symbols.

\smallskip

The \emph{(Hamming) distance}\index{Hamming distance} between two
codewords $\mathbf{x},\mathbf{y} \in \mathbb{F}^n$ is defined by the
number of coordinate positions in which they differ, i.e.
\[d(\mathbf{x},\mathbf{y}):= \left|\{i \mid 1 \leq i \leq n, \;\, x_i \neq y_i\}\right|.\]
The \emph{weight}\index{weight} $w(\mathbf{x})$ of a codeword
$\mathbf{x}$ is defined by
\[w(\mathbf{x}):=d(\mathbf{x},\mathbf{0}),\]
whenever $0$ is an element of $\mathbb{F}$.

\smallskip

A subset $C \subseteq \mathbb{F}^n$ is called a (\emph{$q$-ary})
\emph{code} of \emph{length} $n$ (\emph{binary} if $q=2$,
\emph{ternary} if $q=3$). The elements of $C$ are called
\emph{codewords}. A \emph{linear code}\index{linear code} (or
\emph{$[n,k]$ code}\index{$[n,k]$ code}) over the field
$\mathbb{F}_q$ is a \mbox{$k$-dimensional} linear subspace $C$ of
the vector space $\mathbb{F}_q^n$. We note that large parts of
coding theory are concerned with linear codes. In particular, as
many combinatorial configurations can be described by their
incidence matrices, coding theorists have started in the early
1960's to consider as codes the vector spaces spanned by the rows of
the respective incidence matrices over some given field.

\pagebreak

The \emph{minimum distance}\index{minimum distance} $d$ of a code
$C$ is defined as
\[d:=\;\mbox{min}\,\{d(\mathbf{x},\mathbf{y}) \mid \mathbf{x},\mathbf{y} \in C,\;\,\mathbf{x} \neq \mathbf{y}\}.\]
Clearly, the minimum distance of a linear code is equal to its
\emph{minimum weight}\index{minimum weight}, i.e. the minimum of the
weights of all non-zero codewords. An $[n,k,d]$ code is an $[n,k]$
code with minimum distance $d$.

\smallskip

The minimum distance of a (not necessarily linear) code $C$
determines the error-correcting capability of $C$: If $d=2e+1$, then
$C$ is called an \emph{$e$-error-correcting
code}\index{error-correcting code}. Defining by
\[S_e(\mathbf{x}):=\{\mathbf{y} \in
\mathbb{F}^n \mid d(\mathbf{x},\mathbf{y}) \leq e\}\] the
\emph{sphere}\index{sphere} (or \emph{ball}\index{ball}) of radius
$e$ around a codeword $\mathbf{x}$ of $C$, this implies that the
spheres of radius $e$ around distinct codewords are disjoint.

\smallskip

Counting the number of codewords in a sphere of radius $e$ yields to
the subsequent \emph{sphere packing}\index{sphere packing bound} (or
\emph{Hamming}) \emph{Bound}\index{Hamming Bound}.

\begin{theorem}\label{Hamming-bound}
Let $C$ be a $q$-ary code of length $n$ and minimum distance
$d=2e+1$. Then
\[\left| C \right| \cdot \sum_{i=0}^{e} {n \choose i} (q-1)^i \leq
q^n.\]
\end{theorem}

If equality holds, then $C$ is called a \emph{perfect
code}\index{perfect code}. Equivalently, $C$ is perfect if the
spheres of radius $e$ around all codewords cover the whole space
$\mathbb{F}^n$. Certainly, perfect codes are combinatorially
interesting objects, however, they are extremely rare.

\smallskip

We will call two codes \emph{(permutation)
equivalent}\index{permutation equivalent} if one is obtained from
the other by applying a fixed permutation to the coordinate
positions for all codewords. A \emph{generator
matrix}\index{generator matrix} $G$ for an $[n,k]$ code $C$ is a $(k
\times n)$-matrix for which the rows are a basis of $C$. We say that
$G$ is in \emph{standard form}\index{standard form} if $G=(I_k, P)$,
where $I_k$ is the $(k \times k)$ identity matrix.

\smallskip

For an $[n,k]$ code $C$, let \[C^\perp :=\{\textbf{x} \in
\mathbb{F}_q^n \mid \forall_{\textbf{y} \in C} [\langle
\textbf{x},\textbf{y}\rangle=0]\}\] denote the \emph{dual
code}\index{dual code} of $C$, where $\langle
\textbf{x},\textbf{y}\rangle$ is the standard inner (or dot) product
in $\mathbb{F}_q^n$. The code $C^\perp$ is an $[n,n-k]$ code. If $H$
is a generator matrix for $C^\perp$, then clearly
\[C=\{\textbf{x} \in \mathbb{F}_q^n \mid \textbf{x}H^T = \textbf{0}\},\]
and $H$ is called a \emph{parity check matrix}\index{parity check
matrix} for the code $C$. If $G=(I_k, P)$ is a generator matrix of
$C$, then $H=(-P^T, I_{n-k})$ is a parity check matrix of $C$. A
code $C$ is called \emph{self-dual}\index{self-dual code} if
$C=C^\perp$. If $C \subset C^\perp$, then $C$ is called
\emph{self-orthogonal}\index{self-orthogonal code}.

\smallskip

If $C$ is a linear code of length $n$ over $\mathbb{F}_q$, then
\[\overline{C}:=\{(c_1, \ldots, c_n,c_{n+1}) \mid (c_1, \ldots, c_n) \in C,\; \sum_{i=1}^{n+1} c_i =0\}\]
defines the \emph{extended code}\index{extended code} corresponding
to $C$. The symbol $c_{n+1}$ is called the \emph{overall parity
check symbol}\index{overall parity check symbol}. Conversely, $C$ is
the \emph{punctured} (or \emph{shortened})\index{punctured code}
code of $\overline{C}$.

\smallskip

The \emph{weight distribution}\index{weight distribution} of a
linear code $C$ of length $n$ is the sequence $\{A_i\}_{i=0}^{n}$,
where $A_i$ denotes the number of codewords in $C$ of weight $i$.
The polynomial \[A(x) := \sum_{i=0}^{n} A_i x^i\] is the
\emph{weight enumerator}\index{weight enumerator} of $C$.

\smallskip

The weight enumerators of a liner code and its dual code are
related, as shown by the following theorem, which is one of the most
important results in the theory of error-correcting codes.

\begin{theorem}{\em (MacWilliams~\cite{will63}).}\label{MacWill}
Let $C$ be an $[n,k]$ code over $\mathbb{F}_q$ with weight
enumerator $A(x)$ and let $A^\perp (x)$ be the weight enumerator of
the dual code $C^\perp$. Then
\[A^\perp(x) = q^{-k}(1+(q-1)x)^n A \big(\frac{1-x}{1+(q-1)x} \big).\]
\end{theorem}

We note that the concept of the weight enumerator can be generalized
to non-linear codes (so-called \emph{distance enumerator},
cf.~\cite{will72,dels73} and Subsection~\ref{schemes}).

\smallskip

An $[n,k]$ code $C$ over $\mathbb{F}_q$ is called
\emph{cyclic}\index{cyclic code} if
\[\forall_{(c_0,c_1,\ldots,c_{n-1}) \in C} [(c_{n-1},c_0,\ldots,c_{n-2}) \in C],\]
i.e. any cyclic shift of a codeword is again a codeword. We adopt
the usual convention for cyclic codes that $n$ and $q$ are coprime.
Using the isomorphism
\[(a_0,a_1,\ldots,a_{n-1}) \rightleftarrows a_0 + a_1x + \ldots + a_{n-1}x^{n-1}\]
between $\mathbb{F}_q^n$ and the residue class ring $\mathbb{F}_q[x]
/(x^n-1)$, it follows that a cyclic code corresponds  to an ideal in
$\mathbb{F}_q[x] /(x^n-1)$.

%%% ----------------------------------------------------------------------

\section{Thoughts for Practitioners}\label{pract}

In the following, we introduce and elaborate on the fruitful
interplay of coding theory and algebraic combinatorics, with most of
the focus on the interaction of codes with combinatorial designs,
finite geometries, simple groups, sphere packings, kissing numbers,
lattices, and association schemes. In particular, special interest
is devoted to the relationship between codes and combinatorial
designs. We give an accessible and reasonably self-contained
exposition in the first subsection as we do not assume the reader is
familiar with the theory of combinatorial designs. In what follows,
we describe and recapitulate important results in the development of
the state of the art. In addition, we give illustrative examples and
constructions, and highlight recent achievements.

%%% ----------------------------------------------------------------------

\subsection{Introduction to finite projective planes and combinatorial designs}\label{intro_des}

Combinatorial design theory\index{combinatorial design theory} is a
subject of considerable interest in discrete mathematics. We give in
this subsection an introduction to the topic, with emphasis on the
construction of some important designs. For a more general treatment
of combinatorial designs, the reader is referred
to~\cite{bjl99,cam76,crc06,hall86,hupi85,stin04}. In
particular,~\cite{bjl99,crc06} provide encyclopedias on key results.

\smallskip

Besides coding theory\index{coding theory}, there are many
interesting connections of design theory to other fields. We mention
in our context especially its links to finite
geometries~\cite{demb68}\index{finite geometries}, incidence
geometry~\cite{buek_handb95}\index{incidence geometry}, group
theory~\cite{cam99a,carm37,dimo96,wiel64}\index{group theory}, graph
theory~\cite{cali91,ton88}\index{graph theory},
cryptography~\cite{pei06,stin93,stin99}\index{cryptography}, as well
as classification algorithms~\cite{kaos06}\index{classification
algorithm}. In addition to that, we
recommend~\cite{cam94,wils01,hall86,rys63,grah95} for the reader
interested in the broad area of combinatorics in general.

\smallskip

We start by introducing several notions.

\begin{definition}\label{projplane}\index{projective plane}
A \emph{projective plane of order $n$} is a pair of \emph{points}
and \emph{lines} such that the following properties hold:
\begin{enumerate}
\item[(i)] any two distinct points are on a unique line,

\item[(ii)] any two distinct lines intersect in a unique point,

\item[(iii)] there exists a \emph{quadrangle}\index{quadrangle}, i.e. four points no three of which are on a common line,

\item[(iv)] there are $n+1$ points on each line, $n+1$ lines through
each point and the total number of points, respectively lines, is
$n^2+n+1$.
\end{enumerate}
\end{definition}
It follows easily from (i), (ii), and (iii) that the number of
points on a line is a constant. When setting this constant equal to
$n+1$, then (iv) is a consequence of (i) and (iii).

\smallskip

Combinatorial designs can be regarded as generalizations of
projective planes:

\begin{definition}\label{StDes}
For positive integers $t \leq k \leq v$ and $\lambda$, we define a
\mbox{\emph{$t$-design}}\index{design}\index{$t$-design}, or more
precisely a \emph{$t$-$(v,k,\lambda)$ design}, to be a pair
\mbox{$\mathcal{D}=(X,\mathcal{B})$}, where $X$ is a finite set of
\emph{points}\index{point}, and $\mathcal{B}$ a set of
\mbox{$k$-element} subsets of $X$ called \emph{blocks}\index{block},
with the property that any $t$ points are contained in precisely
$\lambda$ blocks.
\end{definition}
We will denote points by lower-case and blocks by upper-case Latin
letters. Via convention, we set $v:=\left| X \right|$ and $b:=
\left| \mathcal{B} \right|$. Throughout this chapter, `repeated
blocks' are not allowed, that is, the same \mbox{$k$-element} subset
of points may not occur twice as a block. If $t<k<v$ holds, then we
speak of a \emph{non-trivial}\index{non-trivial
design}\index{trivial design} \mbox{$t$-design}.

\smallskip

Designs may be represented algebraically in terms of incidence
matrices: Let $\mathcal{D}=(X,\mathcal{B})$ be a \mbox{$t$-design},
and let the points be labeled $\{x_1,\ldots,x_v\}$ and the blocks be
labeled $\{B_1,\ldots,B_b\}$. Then, the $(b \times v)$-matrix
$A=(a_{ij})$ ($1\leq i \leq b$, $1\leq j \leq v$) defined by
\[a_{ij}:= \left\{\begin{array}{ll}
    1,\;\mbox{if}\; x_j \in B_i \\
    0,\;\mbox{otherwise}\\
\end{array} \right.\]
is called an \emph{incidence matrix}\index{incidence matrix} of
$\mathcal{D}$. Clearly, $A$ depends on the respective labeling,
however, it is unique up to column and row permutation.

\smallskip

If $\mathcal{D}_1=(X_1,\mathcal{B}_1)$ and $\mathcal{D}_2=(X_2,
\mathcal{B}_2)$ are two \mbox{$t$-designs}, then a bijective map
$\alpha: X_1 \longrightarrow X_2$ is called an
\emph{isomorphism}\index{isomorphism} of $\mathcal{D}_1$ onto
$\mathcal{D}_2$, if
\[B \in \mathcal{B}_1 \iff \alpha(B) \in \mathcal{B}_2.\]
In this case, the designs $\mathcal{D}_1$ and $\mathcal{D}_2$ are
\emph{isomorphic}. An isomorphism of a design $\mathcal{D}$ onto
itself is called an \emph{automorphism}\index{automorphism} of
$\mathcal{D}$. Evidently, the set of all automorphisms of a design
$\mathcal{D}$ form a group under composition, the \emph{full
automorphism group}\index{full automorphism group} of $\mathcal{D}$.
Any subgroup of it will be called \emph{an automorphism
group}\index{automorphism group} of $\mathcal{D}$.

\smallskip

If $\mathcal{D}=(X,\mathcal{B})$ is a \mbox{$t$-$(v,k,\lambda)$}
design with $t \geq 2$, and $x \in X$ arbitrary, then the
\emph{derived design with respect to $x$}\index{derived design} is
\mbox{$\mathcal{D}_x=(X_x,\mathcal{B}_x)$}, where $X_x = X
\backslash \{x\}$, \mbox{$\mathcal{B}_x=\{B \backslash \{x\} \mid x
\in B \in \mathcal{B}\}$}. In this case, $\mathcal{D}$ is also
called an \emph{extension}\index{extension} of $\mathcal{D}_x$.
Obviously, $\mathcal{D}_x$ is a \mbox{$(t-1)$-$(v-1,k-1,\lambda)$}
design. The \emph{complementary design} $\mathcal{\overline{D}}$ is
obtained by replacing each block of $\mathcal{D}$ by its complement.

\smallskip

For historical reasons, a \mbox{$t$-$(v,k,\lambda)$ design} with
$\lambda =1$ is called a \emph{Steiner
\mbox{$t$-design}}\index{Steiner design}\index{Steiner $t$-design}.
Sometimes this is also known as a \emph{Steiner
system}\index{Steiner system} if the parameter $t$ is clearly given
from the context.

\smallskip

The special case of a Steiner design with parameters $t=2$ and $k=3$
is called a \emph{Steiner triple system of order $v$} (briefly
$STS(v)$)\index{Steiner triple system}. The question regarding their
existence was posed in the classical ``Combinatorische Aufgabe''
(1853) of the nineteenth century geometer Jakob
Steiner~\cite{stein53}:

\begin{quote}

``Welche Zahl, $N$, von Elementen hat die Eigenschaft, dass sich die
Elemente so zu dreien ordnen lassen, dass je zwei in einer, aber nur
in einer Verbindung vorkommen?''

\end{quote}

However, there had been earlier work on these particular designs
going back to, in particular, J.~Pl\"{u}cker, W.~S.~B.~Woolhouse,
and most notably T.~P.~Kirkman. For an account on the early history
of designs, see~\cite[Chap.\,I.2]{crc06} and~\cite{wils03}.

\smallskip

A Steiner design with parameters $t=3$ and $k=4$ is called a
\emph{Steiner quadruple system of order $v$} (briefly
$SQS(v)$).\index{Steiner quadruple system}

\smallskip

If a \mbox{$2$-design} has equally many points and blocks, i.e.
$v=b$, then we speak of a \emph{square design}\index{square design}
(as its incidence matrix is square). By tradition, square designs
are often called \emph{symmetric designs},\index{symmetric design}
although here the term does not imply any symmetry of the design.
For more on these interesting designs, see, e.g.,~\cite{iosh06}.

\pagebreak

We give some illustrative examples of finite projective planes and
combinatorial designs. We assume that $q$ is always a prime power.

\begin{example}\label{STS7}\index{Steiner triple system}
Let us choose as point set \[X=\{1,2,3,4,5,6,7\}\] and as block set
\[\mathcal{B}=\{\{1,2,4\},\{2,3,5\},\{3,4,6\},\{4,5,7\},\{1,5,6\},\{2,6,7\},\{1,3,7\}\}.\]
This gives a \mbox{$2$-$(7,3,1)$ design}, the well-known \emph{Fano
plane}\index{Fano plane}, the smallest design arising from a
projective geometry, which is unique up to isomorphism. We give the
usual representation of this projective plane of order $2$ by the
following diagram:
\end{example}

\begin{figure}[htp]

\centering

\begin{tikzpicture}[scale=1,thick]

\filldraw [draw=black!100,fill=black!0]

(0,0.87) circle (24.4pt);

\filldraw [draw=black!100,fill=black!100]

(0,0) circle (2pt)  (0,-0.3) node {7}

(-1.5,0) circle (2pt) (-1.52,-0.3) node {1}

(1.5,0) circle (2pt) (1.52,-0.3) node {3}

(0,0.87) circle (2pt) (0.12,1.15) node {6}

(0,2.6) circle(2pt) (0,2.9) node {2}

(0.75,1.31) circle (2pt) (0.95,1.35) node {5}

(-0.75,1.31) circle(2pt) (-0.95,1.35) node {4};

\draw (1.5,0) -- (-1.5,0)

(1.5,0) -- (0,2.6)

(-1.5,0) -- (0,2.6)

(-1.5,0) -- (0.75,1.31)

(1.5,0) -- (-0.75,1.31)

(0,0) -- (0,2.6);

\end{tikzpicture}
\caption{Fano plane}\label{fano}\index{Fano plane}
\end{figure}
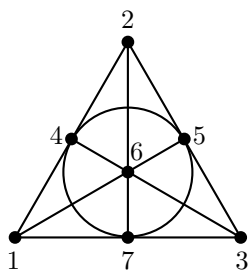

\begin{figure}[htp]
\centering
\begin{tikzpicture}[scale=0.85,thick]

\filldraw [draw=black!100,fill=black!100]

(0,0) circle (2pt) (0,-0.3) node {1}

(0,1) circle (2pt) (-0.25,1) node {4}

(0,2) circle(2pt) (0,2.3) node {7}

(1,0) circle (2pt) (1,-0.3) node {2}

(1,1) circle (2pt) (1.3,1.12) node {5}

(1,2) circle(2pt) (1,2.3) node {8}

(2,0) circle (2pt) (2,-0.3) node {3}

(2,1) circle(2pt) (2.25,1) node {6}

(2,2) circle (2pt) (2,2.3) node {9};

\draw (0,0) -- (0,2)

(0,0) -- (2,0)

(2,0) -- (2,2)

 (0,2) -- (2,2)

 (0,1) -- (2,1)

(1,0) -- (1,2)

(0,0) -- (2,2)

(2,0) -- (0,2)

(0,1) -- (1,0)

(2,1) -- (1,2)

(0,1) -- (1,2)

(1,0) -- (2,1)

(2,1) .. controls (4,3) and (2,4) .. (0,2)

(1,2) .. controls (-1,4) and (-2,2) .. (0,0)

(0,1) .. controls (-2,-1) and (0,-2) .. (2,0)

(1,0) .. controls (3,-2) and (4,0) .. (2,2);

\end{tikzpicture}
\caption{Affine plane of order $3$}\label{AG23}
\end{figure}
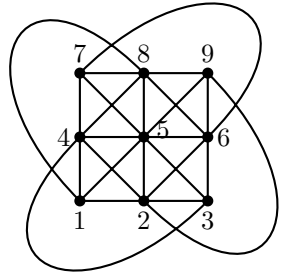

\begin{example}\label{STS9}
We take as point set \[X=\{1,2,3,4,5,6,7,8,9\}\] and as block set
\[\mathcal{B}=\{\{1,2,3\},\{4,5,6\},\{7,8,9\},\{1,4,7\},\{2,5,8\},\{ 3,6,9\},\]
\[ \qquad \quad \{1,5,9\},\{2,6,7\},\{3,4,8\},\{1,6,8\},\{2,4,9\},\{3,5,7\}\}.\]
This gives a \mbox{$2$-$(9,3,1)$ design}, the smallest non-trivial
design arising from an affine geometry, which is again unique up to
isomorphism. This affine plane of order $3$ can be constructed from
the array
\[\begin{array}{ccc}
  1 & 2 & 3 \\
  4 & 5 & 6 \\
  7 & 8 & 9
\end{array}\]
as shown in Figure~\ref{AG23}.
\end{example}

More generally, we obtain:

\begin{example}\index{projective geometry}\label{PG}
We choose as point set $X$ the set of \mbox{$1$-dimensional}
subspaces of a vector space $V=V(d,q)$ of dimension $d \geq 3$ over
$\mathbb{F}_q$. As block set $\mathcal{B}$ we take the set of
\mbox{$2$-dimensional} subspaces of $V$. Then there are $v=(q^{d}
-1)/(q-1)$ points and each block $B \in \mathcal{B}$ contains
$k=q+1$ points. Since obviously any two \mbox{$1$-dimensional}
subspaces span a single \mbox{$2$-dimensional} subspace, any two
distinct points are contained in a unique block. Thus, the
\emph{projective space} $PG(d-1,q)$\index{projective space} is an
example of a \mbox{$2$-$(\frac{q^{d}-1}{q-1},q+1,1)$ design}. For
$d=3$, the particular designs are \emph{projective planes of order
$q$}\index{projective plane}, which are square designs. More
generally, for any fixed $i$ with $1 \leq i \leq d-2$, the points
and \mbox{$i$-dimensional} subspaces of $PG(d-1,q)$ (i.e. the
\mbox{$(i+1)$-dimensional} subspaces of $V$) yield a
\mbox{$2$-design}.
\end{example}

\begin{example}\index{affine geometry}\label{AG}
We take as point set $X$ the set of elements of a vector space
$V=V(d,q)$ of dimension $d \geq 2$ over $\mathbb{F}_q$. As block set
$\mathcal{B}$ we choose the set of affine lines of $V$ (i.e. the
translates of \mbox{$1$-dimensional} subspaces of $V$). Then there
are $v=q^d$ points and each block $B \in \mathcal{B}$ contains $k=q$
points. As clearly any two distinct points lie on exactly one line,
they are contained in a unique block. Hence, we obtain the
\emph{affine space} $AG(d,q)$\index{affine space} as an example of a
\mbox{$2$-$(q^d,q,1)$ design}. When $d=2$, these designs are
\emph{affine planes of \mbox{order $q$}}\index{affine plane}. More
generally, for any fixed $i$ with $1 \leq i \leq d-1$, the points
and \mbox{$i$-dimensional} subspaces of $AG(d,q)$ form a
\mbox{$2$-design}.
\end{example}

\begin{remark}\label{projpl_existence}
It is well-established that both affine and projective planes of
order $n$ exist whenever $n$ is a prime power. The conjecture that
no such planes exist with orders other than prime powers is
unresolved so far. The classical result of R.~H.~Bruck and
H.~J.~Ryser~\cite{brry49} still remains the only general statement:
If $n \equiv 1$ or $2$ (mod $4)$ and $n$ is not equal to the sum of
two squares of integers, then $n$ does not occur as the order of a
finite projective plane. The smallest integer that is not a prime
power and not covered by the Bruck-Ryser Theorem\index{Bruck-Ryser
Theorem} is $10$. Using substantial computer analysis, C.~W.~H.~Lam,
L.~Thiel, and S.~Swiercz~\cite{lam89} proved the non-existence of a
projective plane of order $10$ (cf.~Remark~\ref{projpl10}). The next
smallest number to consider is $12$, for which neither a positive
nor a negative answer has been proved.

Needless to mention that --- apart from the existence problem ---
the question on the number of different isomorphism types (when
existent) is fundamental. There are, for example, precisely four
non-isomorphic projective planes of order $9$. For a further
discussion, in particular of the rich history of affine and
projective planes, we refer, e.g.,
to~\cite{beut95,demb68,hirsch98,hupi82,luene80,pick75}.
\end{remark}

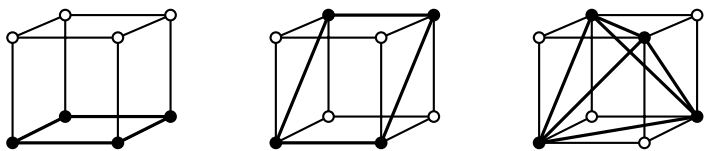
\begin{figure}[htp]
\centering
\begin{tikzpicture}[scale=1,thick]

\filldraw [draw=black!100,fill=black!100]

(0,0) circle (2pt)

(1.4,0) circle (2pt)

(0.7,0.35) circle (2pt)

(2.1,0.35) circle (2pt);

\draw

(0,0) -- (0,1.4)

(0,0) -- (0.7,0.35)

(1.4,0) -- (1.4,1.4)

(1.4,0) -- (2.1,0.35)

(0.7,0.35) -- (2.1,0.35)

(2.1,0.35) -- (2.1,1.7)

(0,1.4) -- (1.4,1.4)

(0.7,1.7) -- (2.1,1.7)

(1.4,1.4) -- (2.1,1.7)

(0,1.4) -- (0.7,1.7)

(0.7,0.35) -- (0.7,1.7);

\filldraw [draw=black!100,fill=black!0]

(0,1.4) circle(2pt)

(1.4,1.4) circle (2pt)

(0.7,1.7) circle(2pt)

(2.1,1.7) circle (2pt);

\draw[very thick]

(0,0) -- (1.4,0)

(0,0) -- (0.7,0.35)

(1.4,0) -- (2.1,0.35)

(0.7,0.35) -- (2.1,0.35);

\filldraw [draw=black!100,fill=black!100]

(3.5,0) circle (2pt)

(4.9,0) circle (2pt)

(4.2,1.7) circle(2pt)

(5.6,1.7) circle (2pt);

\draw (3.5,0) -- (4.9,0)

(3.5,0) -- (3.5,1.4)

(3.5,0) -- (4.2,0.35)

(4.9,0) -- (4.9,1.4)

(4.9,0) -- (5.6,0.35)

(4.2,0.35) -- (5.6,0.35)

(5.6,0.35) -- (5.6,1.7)

(3.5,1.4) -- (4.9,1.4)

(4.2,1.7) -- (5.6,1.7)

(4.9,1.4) -- (5.6,1.7)

(3.5,1.4) -- (4.2,1.7)

(4.2,0.35) -- (4.2,1.7);

\filldraw [draw=black!100,fill=black!0]

(4.2,0.35) circle (2pt)

(5.6,0.35) circle (2pt)

(3.5,1.4) circle(2pt)

(4.9,1.4) circle (2pt);

\draw[very thick]

(3.5,0) -- (4.9,0)

(3.5,0) -- (4.2,1.7)

(4.9,0) -- (5.6,1.7)

(4.2,1.7) -- (5.6,1.7);

\filldraw [draw=black!100,fill=black!100]

(7,0) circle (2pt)

(7.7,1.7) circle(2pt)

(8.4,1.4) circle (2pt)

(9.1,0.35) circle (2pt);

\draw (7,0) -- (8.4,0)

(7,0) -- (7,1.4)

(7,0) -- (7.7,0.35)

(8.4,0) -- (8.4,1.4)

(8.4,0) -- (9.1,0.35)

(7.7,0.35) -- (9.1,0.35)

(9.1,0.35) -- (9.1,1.7)

(7,1.4) -- (8.4,1.4)

(7.7,1.7) -- (9.1,1.7)

(8.4,1.4) -- (9.1,1.7)

(7,1.4) -- (7.7,1.7)

(7.7,0.35) -- (7.7,1.7);

\filldraw [draw=black!100,fill=black!00]

(8.4,0) circle (2pt)

(7.7,0.35) circle (2pt)

(7,1.4) circle(2pt)

(9.1,1.7) circle (2pt);

\draw[very thick]

(7,0) -- (7.7,1.7)

(7,0) -- (8.4,1.4)

(7,0) -- (9.1,0.35)

(7.7,1.7) -- (9.1,0.35)

(8.4,1.4) -- (9.1,0.35)

(8.4,1.4) -- (7.7,1.7);
\end{tikzpicture}

\caption{Steiner quadruple system of order $8$}\label{cube}

\end{figure}

\begin{example}\label{AG_cube}
We take as points the vertices of a \mbox{$3$-dimensional} cube. As
illustrated in Figure~\ref{cube}\index{cube}, we can choose three
types of blocks:
\begin{enumerate}
\item[(i)] a face (six of these),
\item[(ii)] two opposite edges (six of these),
\item[(iii)] an inscribed regular tetrahedron\index{regular tetrahedron} (two of these).
\end{enumerate}
This gives a \mbox{$3$-$(8,4,1)$ design}, which is unique up to
isomorphism.
\end{example}

We have more generally:

\begin{example}\label{AGd2}
In $AG(d,q)$ any three distinct points define a plane unless they
are collinear\index{collinear} (that is, lie on the same line). If
the underlying field is $\mathbb{F}_2$, then the lines contain only
two points and hence any three points cannot be collinear.
Therefore, the points and planes of the affine space $AG(d,2)$ form
a \mbox{$3$-$(2^d,4,1)$} design. More generally, for any fixed $i$
with $2 \leq i \leq d-1$, the points and \mbox{$i$-dimensional}
subspaces of $AG(d,2)$ form a \mbox{$3$-design}.
\end{example}

\begin{example}\label{Witt_des}\index{Mathieu-Witt designs}
The unique \mbox{$2$-$(9,3,1)$ design} whose points and blocks are
the points and lines of the affine plane $AG(2,3)$ can be extended
precisely three times to the following designs which are also unique
up to isomorphism: the \mbox{$3$-$(10,4,1)$ design} which is the
M{\"{o}}bius plane of order $3$ with $P \mathit{\Gamma} L(2,9)$ as
full automorphism group, and the two \emph{Mathieu-Witt designs}
\mbox{$4$-$(11,5,1)$} and \mbox{$5$-$(12,6,1)$} with the sporadic
Mathieu groups\index{Mathieu groups} $M_{11}$ and $M_{12}$ as point
\mbox{$4$-transitive} and point \mbox{$5$-transitive} full
automorphism groups, respectively.

\pagebreak

To construct the `large' Mathieu-Witt designs one starts with the
\mbox{$2$-$(21,5,1)$ design} whose points and blocks are the points
and lines of the projective plane $PG(2,4)$. This can be extended
also exactly three times to the following unique designs: the
\emph{Mathieu-Witt design} \mbox{$3$-$(22,6,1)$} with
$\mbox{Aut}(M_{22})$ as point \mbox{$3$-transitive} full
automorphism group as well as the \emph{Mathieu-Witt designs}
\mbox{$4$-$(23,7,1)$} and \mbox{$5$-$(24,8,1)$} with $M_{23}$ and
$M_{24}$ as point \mbox{$4$-transitive} and point
\mbox{$5$-transitive} full automorphism groups, respectively.

The five Mathieu groups were the first sporadic simple
groups\index{sporadic simple groups} and were discovered by
E.~Mathieu~\cite{math61,math73} over one hundred years ago. They are
the only finite \mbox{$4$-} and \mbox{$5$-transitive} permutation
groups apart from the symmetric or alternating groups. The Steiner
designs associated with the Mathieu groups were first constructed by
both R.~D.~Carmichael~\cite{carm37} and E.~Witt~\cite{witt38a}, and
their uniqueness established up to isomorphism by
Witt~\cite{witt38}. From the meanwhile various alternative
constructions, we mention especially those of
H.~L\"{u}neburg~\cite{luene69} and
M.~Aschbacher~\cite[Chap.\,6]{asch94}. However, the easiest way to
construct and prove uniqueness of the Mathieu-Witt designs is via
coding theory, using the related \emph{binary} and \emph{ternary}
\emph{Golay codes} (see Subsection~\ref{golay_witt_mathieu}).
\end{example}

\begin{remark}
By classifying Steiner designs which admit automorphism groups with
sufficiently strong symmetry properties, specific characterizations
of the Mathieu-Witt designs with their related Mathieu groups were
obtained (see, e.g.,~\cite{tits64,luene65,hu01,hu05,hu07a,hu07}
and~\cite[Chap.\,5]{hu08} for a survey).
\end{remark}

\begin{remark}
We mention that, in general, for $t=2$ and $3$, there are many
infinite classes of Steiner \mbox{$t$-designs}, but for $t=4$ and
$5$ only a finite number are known. Although
L.~Teirlinck~\cite{teir87} has shown that non-trivial
\mbox{$t$-designs} exist for all values of $t$, no Steiner
\mbox{$t$-designs} have been constructed for $t \geq 6$ so far.
\end{remark}

\smallskip

In what follows, we need some helpful combinatorial tools:

\smallskip

A standard combinatorial double counting argument gives the
following assertions.
\begin{lemma}\label{s-design}
Let $\mathcal{D}=(X,\mathcal{B})$ be a \mbox{$t$-$(v,k,\lambda)$}
design. For a positive integer $s \leq t$, let $S \subseteq X$ with
$\left|S\right|=s$. Then the total number $\lambda_s$ of blocks
containing all the points of $S$ is given by
\[\lambda_s = \lambda \frac{{v-s \choose t-s}}{{k-s \choose t-s}}.\]
In particular, for $t\geq 2$, a \mbox{$t$-$(v,k,\lambda)$} design is
also an \mbox{$s$-$(v,k,\lambda_s)$} design.
\end{lemma}

For historical reasons, it is customary to set $r:= \lambda_1$ to be
the total number of blocks containing a given point (referring to
the `replication number' from statistical design of
experiments\index{statistical design of experiments}, one of the
origins of designs theory).

\pagebreak

\begin{lemma}\label{Comb_t=5}
Let $\mathcal{D}=(X,\mathcal{B})$ be a \mbox{$t$-$(v,k,\lambda)$}
design. Then the following holds:
\begin{enumerate}

\item[\emph{(a)}] $bk = vr.$

\smallskip

\item[\emph{(b)}] $\displaystyle{{v \choose t} \lambda = b {k \choose t}.}$

\smallskip

\item[\emph{(c)}] $r(k-1)=\lambda_2(v-1)$ for $t \geq 2$.

\end{enumerate}
\end{lemma}

Since in Lemma~\ref{s-design} each $\lambda_s$ must be an integer,
we have moreover the subsequent necessary arithmetic conditions.

\begin{lemma}\label{divCond}
Let $\mathcal{D}=(X,\mathcal{B})$ be a \mbox{$t$-$(v,k,\lambda)$}
design. Then
\[\lambda {v-s \choose t-s} \equiv \, 0\; \emph{(mod}\;\, {k-s \choose t-s})\]
for each positive integer $s \leq t$.
\end{lemma}

The following theorem is an important result in the theory of
designs, generally known as \emph{Fisher's
Inequality}\index{Fisher's Inequality}.

\begin{theorem}{\em (Fisher~\cite{fish40}).}\label{FisherIn}
If $\mathcal{D}=(X,\mathcal{B})$ is a non-trivial
\mbox{$t$-$(v,k,\lambda)$} design with $t \geq 2$, then we have $b
\geq v$, that is, there are at least as many blocks as points in
$\mathcal{D}$.
\end{theorem}

We remark that equality holds exactly for square designs when $t=2$.
Obviously, the equality $v=b$ implies $r=k$ by
Lemma~\ref{Comb_t=5}~(a).

%%% ----------------------------------------------------------------------

\subsection{Basic connections between codes and combinatorial designs}\label{basic}

There is a rich and fruitful interplay between coding theory and
design theory. In particular, many \mbox{$t$-designs} have been
found in the last decades by considering the codewords of fixed
weight in some special, often linear codes. As we will see in the
sequel, these codes typically exhibit a high degree of regularity.
There is an amount of
literature~\cite{assma74,assk93,assk96,bla79,cali91,hupl03,lint77,lint93,ton88,ton98,ton06,wisl77}
discussing to some extent in more detail various relations between
codes and designs.

\smallskip

For a codeword $x \in \mathbb{F}^n$, the set
\[\mbox{supp}(x):=\{i \mid x_i \neq 0\}\] of all coordinate positions with non-zero
coordinates is called the \emph{support}\index{support} of $x$. We
shall often form a \mbox{$t$-design} of a code in the following way:
Given a (usually linear) code of length $n$, which contains the zero
vector, and non-zero weight $w$, we choose as point set $X$ the set
of $n$ coordinate positions of the code and as block set
$\mathcal{B}$ the supports of all codewords of weight $w$.

Since we do not allow repeated blocks, clearly only distinct
representatives of supports for codewords with the same supports are
taken in the non-binary case.

\smallskip

\pagebreak

We give some elementary examples.

\begin{example}\label{H7} The matrix
\[G=\left(
      \begin{array}{ccccccc}
        1 & 1 & 0 & 1 & 0 & 0 & 0 \\
        0 & 1 & 1 & 0 & 1 & 0 & 0 \\
        0 & 0 & 1 & 1 & 0 & 1 & 0 \\
        0 & 0 & 0 & 1 & 1 & 0 & 1 \\
      \end{array}
    \right)
 \]
is a generator matrix of a binary $[7,4,3]$ Hamming code, which is
the smallest non-trivial Hamming code (see also
Example~\ref{hamming}). This code is a perfect
single-error-correcting code with weight distribution $A_0=A_7=1$,
$A_3=A_4=7$. The seven codewords of weight $3$ are precisely the
seven rows of the incidence matrix \[\left(
       \begin{array}{ccccccc}
         1 & 1 & 0 & 1 & 0 & 0 & 0 \\
         0 & 1 & 1 & 0 & 1 & 0 & 0 \\
         0 & 0 & 1 & 1 & 0 & 1 & 0 \\
         0 & 0 & 0 & 1 & 1 & 0 & 1 \\
         1 & 0 & 0 & 0 & 1 & 1 & 0 \\
         0 & 1 & 0 & 0 & 0 & 1 & 1 \\
         1 & 0 & 1 & 0 & 0 & 0 & 1 \\
       \end{array}
     \right)\]
of the Fano plane $PG(2,2)$ of Fig.~\ref{fano}. The supports of the
seven codewords of weight $4$ yield the complementary
\mbox{$2$-$(7,4,2)$} design, i.e. the biplane of order $2$.
\end{example}

\begin{example}\label{H8}
The matrix $(I_4,J_4-I_4)$, where $J_4$ denotes the $(4 \times 4)$
all-one matrix, generates the extended binary $[8,4,4]$ Hamming
code. This code is self-dual and has weight distribution
$A_0=A_8=1$, $A_4=14$. As any two codewords of weight $4$ have
distance at least $4$, they have at most two $1$'s in common, and
hence no codeword of weight $3$ can appear as a subword of more than
one codeword. On the other hand, there are ${8\choose 3} =56$ words
of weight $3$ and each codeword of weight $4$ has four subwords of
weight $3$. Hence each codeword of weight $3$ is a subword of
exactly one codeword of weight $4$. Therefore, the supports of the
fourteen codewords of weight $4$ form a \mbox{$3$-$(8,4,1)$} design,
which is the unique $SQS(8)$ (cf.~Example~\ref{AG_cube}).
\end{example}

\smallskip

We give also a basic example of a non-linear code constructed from
design theory.

\begin{example}
We take the rows of an incidence matrix of the (unique) Hadamard
\mbox{$2$-$(11,5,2)$} design, and adjoin the all-one codeword. Then,
the twelve codewords have mutual distance $6$, and if we delete a
coordinate, we get a binary non-linear code of length $10$ and
minimum distance $5$.
\end{example}

For a detailed description of the connection between non-linear
codes and design theory as well as the application of design theory
in the area of (majority-logic) decoding, the reader is referred,
e.g., to~\cite{wisl77,ton98,ton06}.

\smallskip

\pagebreak

Using highly transitive permutation groups, a further construction
of designs from codes can be described (see, e.g.,~\cite{ton88}).

\begin{theorem}\label{t-trans}
Let $C$ be a code which admits an automorphism group acting
\mbox{$t$-homogeneously} (in particular, \mbox{$t$-transitively}) on
the set of coordinates. Then the supports of the codewords of any
non-zero weight form a \mbox{$t$-design}.
\end{theorem}

\begin{example}
The $r$-th order \emph{Reed-Muller (RM) code}\index{Reed-Muller (RM)
code} RM$(r,m)$ of length $2^m$ is a binary $[2^m, \sum_{i=0}^{r} {m
\choose i}, 2^{m-r}]$ code with its codewords the value-vectors of
all Boolean functions in $m$ variables of degree at most $r$. These
codes were first considered by D.~E.~Muller~\cite{mu54} and
I.~S.~Reed~\cite{ree54} in 1954. The dual of the Reed-Muller code
RM$(r,m)$ is RM$(m-r-1,m)$. Clearly, the extended binary $[8,4,4]$
Hamming code in Example~\ref{H8} is RM$(1,3)$.

Alternatively, a codeword in RM$(r,m)$ can be viewed as the sum of
characteristic functions of subspaces of dimension at least $m-r$ of
the affine space $AG(m,2)$. Thus, the full automorphism group of
RM$(r,m)$ contains the \mbox{$3$-transitive} group $AGL(m,2)$ of all
affine transformations, and hence the codewords of any fixed
non-zero weight yield a \mbox{$3$-design}.
\end{example}

%%% ----------------------------------------------------------------------

\subsection{Perfect codes and designs}\label{perf}

The interplay between coding theory and combinatorial designs is
most evidently seen in the relationship between perfect codes and
\mbox{$t$-designs}.

\begin{theorem}{\em (Assmus and Mattson~\cite{assma67}).}\label{assmatt1}
A linear \mbox{$e$-error-correcting} code of length $n$ over
$\mathbb{F}_q$ is perfect if and only if the supports of the
codewords of minimum weight $d=2e+1$ form an
\mbox{$(e+1)$-$(n,d,(q-1)^e)$} design.
\end{theorem}

The question
\begin{quote}
``Does every Steiner triple system on $n$ points extend to a Steiner
quadruple system on $n+1$ points?''
\end{quote}
which goes also back to Jakob Steiner~\cite{stein53}, is still
unresolved in general. However, in terms of binary
$e$-error-correcting codes, there is a positive answer.

\begin{theorem}{\em (Assmus and Mattson~\cite{assma67}).}\label{assmatt2}
Let $C$ be a (not necessarily linear) binary  \mbox{$e$-error
correcting} code of length $n$, which contains the zero vector. Then
$C$ is perfect if and only if the supports of the codewords of
minimum weight $d=2e+1$ form a Steiner \mbox{$(e+1)$-$(n,d,1)$}
design. Moreover, the supports of the minimum codewords in the
extended code $\overline{C}$ form a Steiner
\mbox{$(e+2)$-$(n+1,d+1,1)$} design.
\end{theorem}

We have seen in Example~\ref{H7} and Example~\ref{H8} that the
supports of the seven codewords of weight $3$ in the binary
$[7,4,3]$ Hamming code form a $STS(7)$, while the supports of the
fourteen codewords of weight $4$ in the extended $[8,4,4]$ Hamming
code yield a $SQS(8)$. In view of the above theorems, we get more
generally:

\begin{example}\label{hamming}
Let $n:=(q^m-1)/(q-1)$. We consider a $(m \times n)$-matrix $H$ over
$\mathbb{F}_q$ such that no two columns of $H$ are linearly
dependent. Then $H$ clearly is a parity check matrix of an
$[n,n-m,3]$ code, which is the \emph{Hamming code}\index{Hamming
code} over $\mathbb{F}_q$. The number of its codewords is $q^{n-m}$,
and for any codeword $\mathbf{x}$, we have
$S_1(\mathbf{x})=1+n(q-1)=q^m$. Hence, by the Sphere Packing Bound
(Theorem~\ref{Hamming-bound}), this code is perfect, and the
supports of codewords of minimum weight $3$ form a
\mbox{$2$-$(n,3,q-1)$} design. Furthermore, in a binary
$[2^m-1,2^m-1-m,3]$ Hamming code the supports of codewords of weight
$3$ form a $STS(2^m-1)$, and the supports of the codewords of weight
$4$ in the extended code yield a $SQS(2^m)$.
\end{example}

\begin{note}
The Hamming codes were developed by R.~W.~Hamming~\cite{ham50} in
the mid 1940's, who was employed at Bell Laboratories, and addressed
a need for error correction in his work on the primitive computers
of the time. We remark that the extended binary $[2^m,2^m-m-1,4]$
Hamming code is the Reed-Muller code RM$(m-2,m)$.
\end{note}

\begin{example}\label{golay}
The \emph{binary Golay code}\index{binary Golay code} is a
$[23,12,7]$ code, while the \emph{ternary Golay code}\index{ternary
Golay code} is a $[11,6,5]$ code. For both codes, the parameters
imply equality in the Sphere Packing Bound, and hence these codes
are perfect. We will discuss later various constructions of these
some of the most famous codes (see~Example~\ref{golay_QR} and
Construction~\ref{golay_constr}). By the above theorems, the
supports of codewords of minimum weight $7$ in the binary
$[23,12,7]$ Golay code form a Steiner \mbox{$4$-$(23,7,1)$} design,
and the supports of the codewords of weight $8$ in the extended
binary $[24,12,8]$ Golay code give a Steiner \mbox{$5$-$(24,8,1)$}
design. The supports of codewords of minimum weight $5$ in the
ternary $[11,6,5]$ Golay code yield a \mbox{$3$-$(11,5,4)$} design.
It can be shown (e.g., via Theorem~\ref{t-trans}) that this is
indeed a Steiner \mbox{$4$-$(11,5,1)$} design. We will see in
Example~\ref{assmatt_golay} that the supports of the codewords of
weight $6$ in the extended ternary $[12,6,6]$ Golay code give a
Steiner \mbox{$5$-$(12,6,1)$} design; thus the above results are not
best possible.
\end{example}

\begin{note}
The Golay codes were discovered by M.~J.~E.~Golay~\cite{gol49} in
1949 in the process of extending Hamming's construction. They have
numerous practical real-world applications, e.g., the use of the
extended binary Golay code in the Voyager spacecraft program during
the early 1980's or in contemporary standard Automatic Link
Establishment (ALE) in High Frequency (HF) data communication for
Forward Error Correction (FEC).
\end{note}

\begin{remark}
It is easily seen from their construction that the Hamming codes are
unique (up to equivalence). It was shown by V.~Pless~\cite{pless68}
that this is also true for the Golay codes. Moreover, the binary and
ternary Golay codes are the only non-trivial perfect
$e$-error-correcting codes with $e>1$ over any field $\mathbb{F}_q$.
Using integral roots of the Lloyd polynomial, this remarkable fact
was proven by A.~Tiet{\"a}v{\"a}inen~\cite{tiet73} and J. H. van
Lint~\cite{lint71}, and independently by V.~A.~Zinov'ev and
V.~K.~Leont'ev~\cite{zile73}. M.~R.~Best~\cite{best82} and
Y.~Hong~\cite{hong84} extended this result to arbitrary alphabets
for $e>2$. For a thorough account of perfect codes, we refer
to~\cite{lint75} and~\cite[Chap.\,11]{coh97}.

\smallskip

As trivial perfect codes can only form trivial designs, we have (up
to equivalence) a complete list of non-trivial linear perfect codes
with their associated designs:

\begin{InTextTable}
\begin{tabular}{|l|l|c|l|}
  \hline
  {Code}  &\multicolumn{2}{l|}{Code parameters}  & {Design parameters}\\
  \hline\hline

  Hamming code & $[\frac{q^m-1}{q-1},\frac{q^m-1}{q-1}-m,3]$ & $q$ any prime power & $2$-$(\frac{q^m-1}{q-1},3,q-1)$ \\

  binary Golay code & $[23,12,7]$ & $q=2$ & $4$-$(23,7,1)$ \\

  ternary Golay code & $[11,6,5]$ & $q=3$ & $4$-$(11,5,1)$  \\
  \hline
  \end{tabular}
\label{table:perf}
\end{InTextTable}
There are various constructions of non-linear
single-error-correcting perfect codes. For more details, see,
e.g.,~\cite{wisl77,lint95,ton98,ton06,rom06} and references therein.
However, a classification of these codes seems out of reach at
present, although some progress has been made recently, see, for
instance~\cite{avhe04,hehe06,phel05}.
\end{remark}

\begin{remark}
The long-standing question whether every Steiner triple system of
order $2^m-1$ occurs in a perfect code has been answered recently in
the negative. Relying on the classification~\cite{kop06} of the
Steiner quadruple systems of order $16$, it was shown in~\cite{op07}
that the unique anti-Pasch Steiner triple system of order $15$
provides a counterexample.
\end{remark}

\begin{remark}
Due to the close relationship between perfect codes and some of the
most interesting designs, several natural extensions of perfect
codes have been examined in this respect: \emph{Nearly perfect
codes}\index{nearly perfect code}~\cite{gosn72}, and the more
general class of \emph{uniformly packed codes}\index{uniformly
packed codes}~\cite{sema71,goti75}, were studied extensively and
eventually lead to \mbox{$t$-designs}. H.~C.~A.~van
Tilborg~\cite{til76} showed that \mbox{$e$-error} correcting
uniformly packed codes do not exist for $e> 3$, and classified those
for $e \leq 3$. For more details,
see~\cite{cali91,lint99,til76,wisl77}. The concept of \emph{diameter
perfect codes}\index{diameter perfect code}~\cite{ahls01,schet02} is
related particularly to Steiner designs. For further generalizations
of perfect codes, see e.g.,~\cite[Chap.\,11]{coh97}
and~\cite[Chap.\,6]{wisl77}.
\end{remark}

%%% ----------------------------------------------------------------------

\subsection{The Assmus-Mattson Theorem and analogues}\label{assmatt_thm}

We consider in this subsection one of the most fundamental results
in the interplay of coding theory and design theory. We start by
introducing two important classes of codes.

\smallskip

Let $q$ be an odd prime power. We define a function $\chi$ (the
so-called \emph{Legendre-symbol}) on $\mathbb{F}_q$ by \[\chi(x):=
\left\{\begin{array}{lll}
    0,\;\mbox{if}\; x =0\\
    1,\;\mbox{if}\; x \; \mbox{is a non-zero square}\\
    -1,\;\mbox{otherwise}.\\
\end{array} \right.\]
We note that $\chi$ is a character on the multiplicative group of
$\mathbb{F}_q$. Using the elements of $\mathbb{F}_q$ as row and
column labels $a_i$ and $a_j$ $(0 \leq i,j < q)$, respectively, a
matrix $Q=(q_{ij})$ of order $q$ can be defined by
\begin{equation}\label{circ}
q_{ij}:= \chi(a_j-a_i).
\end{equation}
If $q$ is a prime, then $Q$ is a circulant matrix. We call a matrix
\[C_{q+1}:=\left(
  \begin{array}{cccc}
    0 & \; 1 & \;\cdots & 1 \\
    \chi(-1) &  &  &  \\
    \vdots &  & Q &  \\
    \chi(-1) &  &  &  \\
  \end{array}
\right)\] of order $q+1$ a \emph{Paley matrix}\index{Paley matrix}.
These matrices were constructed by R.~A.~Paley in 1933 and are a
specific type of conference matrices, which have their origin in the
application to conference telephone circuits.

\begin{construction}\label{QR}
Let $n$ be an odd prime, and $q$ be a \emph{quadratic residue} (mod
$n$), i.e. $q^{(n-1)/2} \equiv$ $1$ (mod $n$). The \emph{quadratic
residue code}\index{quadratic residue code} (or \emph{QR
code}\index{QR code}) of length $n$ over $\mathbb{F}_q$ is a
$[n,(n+1)/2]$ code with minimum weight $d \geq \sqrt{n}$ (so-called
\emph{Square Root Bound}\index{Square Root Bound}). It can be
generated by the $(0,1)$-circulant matrix of order $n$ with top row
the incidence vector of the non-zero quadratic residues (mod $n$).
These codes are a special class of cyclic codes and were first
constructed by A.~M.~Gleason in 1964. For $n \equiv$ $3$ (mod $4$),
the extended quadratic residue code is self-dual. We note for the
important binary case that $q$ is a quadratic residue (mod $n$) if
and only if \mbox{$n \equiv$ $\pm1$ (mod $8$).}
\end{construction}

\begin{note}
By a theorem of A.~M.~Gleason and E.~Prange, the full automorphism
group of an extended quadratic residue code of length $n$ contains
the group $PSL(2,n)$ of all linear fractional transformations whose
determinant is a non-zero square.
\end{note}

\begin{example}\label{golay_QR}
The binary $[7,4,3]$ Hamming code is a quadratic residue code of
length $7$ over $\mathbb{F}_2$. The binary $[23,12,7]$ Golay code is
a quadratic residue code of length $23$ over $\mathbb{F}_2$, while
the ternary $[11,6,5]$ Golay code is a quadratic residue code of
length $11$ over $\mathbb{F}_3$.
\end{example}

\begin{construction}\label{pless}
For $q \equiv$ $-1$ (mod $6$) a prime power, the \emph{Pless
symmetry code}\index{Pless symmetry code} Sym$_{2(q+1)}$ of
dimension $q+1$ is a ternary $[2(q+1),q+1]$ code with generator
matrix $G_{2(q+1)}:=(I_{q+1},C_{q+1})$, where $C_{q+1}$ is a Paley
matrix. Since $C_{q+1}C^T_{q+1}= -I_{q+1}$ (over $\mathbb{F}_3$) for
$q \equiv$ $-1$ (mod $3$), the code Sym$_{2(q+1)}$ is self-dual.
This infinite family of cyclic codes were introduced by
V.~Pless~\cite{pless70,pless72} in 1972. We note that the first
symmetry code $S_{12}$ is equivalent to the extended $[12,6,6]$
Golay code.
\end{construction}

The celebrated Assmus-Mattson Theorem gives a sufficient condition
for the codewords of constant weight in a linear code to form a
\mbox{$t$-design}.

\pagebreak

\begin{theorem}{\em (Assmus and Mattson~\cite{assma69}).}
Let $C$ be an $[n,k,d]$ code over $\mathbb{F}_q$ and $C^\perp$ be
the $[n,n-k,e]$ dual code. Let $n_0$ be the largest integer such
that $n_0-\frac{n_0+q-2}{q-1} < d$, and define $m_0$ similarly for
the dual code $C^\perp$, whereas if $q=2$, we assume that
$n_0=m_0=n$. For some integer $t$ with $0<t<d$, let us suppose that
there are at most $d-t$ non-zero weights $w$ in $C^\perp$ with $w
\leq n-t$. Then, for any weight $v$ with $d \leq v \leq n_0$, the
supports of codewords of weight $v$ in $C$ form a \mbox{$t$-design}.
Furthermore, for any weight $w$ with $e \leq w \leq
\mbox{\emph{min}}\{n-t,m_0\}$, the support of the codewords $w$ in
$C^\perp$ also form a \mbox{$t$-design}.
\end{theorem}

The proof of the theorem involves a clever use of the MacWilliams
relations (Theorem~\ref{MacWill}). Along with these,
Lemma~\ref{s-design} and the immediate observation that codewords of
weight less than $n_0$ with the same support must be scalar
multiples of each other, form the basis of the proof (for a detailed
proof, see, e.g.,~\cite[Chap.\,14]{cali91}).

\begin{remark}
Until this result by E.~F.~Assmus{, Jr.} and H.~F.~Mattson{, Jr.} in
1969, only very few \mbox{$5$-designs} were known: the Mathieu-Witt
designs \mbox{$5$-$(12,6,1)$} and \mbox{$5$-$(24,8,1)$}, the
\mbox{$5$-$(24,8,48)$} design formed by the codewords of weight $12$
(the \emph{dodecads}\index{dodecad}) in the extended binary Golay
code, as well as \mbox{$5$-$(12,6,2)$} and \mbox{$5$-$(24,8,2)$}
designs which had been found without using coding theory. However,
by using the Assmus-Mattson Theorem, it was possible to find a
number of new \mbox{$5$-designs}. In particular, the theorem is most
useful when the dual code has relatively few non-zero weights.
Nevertheless, it has not been possible to detect \mbox{$t$-designs}
for $t>5$ by the Assmus-Mattson Theorem.
\end{remark}

We illustrate in the following examples some applications of the
theorem.

\begin{example}\label{assmatt_golay}
The extended binary $[24,12,8]$ Golay code is self-dual
(cf.~Construction~\ref{QR}) and has codewords of weight $0,8,12,16$,
and $24$ in view of Theorem~\ref{MacWill}. For $t=5$, we obtain the
Steiner \mbox{$5$-$(24,8,1)$} design as in Example~\ref{golay}. In
the self-dual extended ternary $[12,6,6]$ Golay code all codewords
are divisible by $3$, and hence for $t=5$, the supports of the
codewords of weight $6$ form a Steiner \mbox{$5$-$(12,6,1)$} design.
\end{example}

\begin{example}\label{assmatt_QR}
The extended quadratic residue code of length $48$ over
$\mathbb{F}_2$ is self-dual with minimum distance $12$. By
Theorem~\ref{MacWill}, it has codewords of weight
$0,12,16,20,24,28,32,36$, and $48$. For $t=5$, each of the values
$v=12,16,20,$ or $24$ yields a different \mbox{$5$-design} and its
complementary design.
\end{example}

\begin{example}\label{assmatt_pless}
The Pless symmetry code Sym$_{36}$ of dimension $18$ is self-dual
(cf.~Construction~\ref{pless}) with minimum distance $12$. The
supports of codewords of weight $12,15,18,$ and $21$ yield
\mbox{$5$-designs} together with their complementary designs.
\end{example}

\begin{remark}
We give an overview of the state of knowledge concerning codes over
$\mathbb{F}_q$ with their associated \mbox{$5$-designs} (cf. also
the tables in~\cite[Chap.\,16]{wisl77}, \cite{assma74,ton98,ton06}).
In fact, these codes are all self-dual. Trivial designs as well as
complementary designs are omitted.

\begin{InTextTable}
\begin{tabular}{|l|l|c|l|l|}
  \hline
  \mbox{Code}  &\multicolumn{2}{l|}{Code parameters}  & Design parameters & Ref.\\
  \hline\hline
  Extended cyclic code        & $[18,9,8]$ &  $q=4$ & $5$-$(18,8,6)$ & \cite{will78} \\
                               &            &        & $5$-$(18,10,180)$ & \\
  \hline
  Extended binary Golay code & $[24,12,8]$ & $q=2$ & $5$-$(24,8,1)$ & \cite{pai56} \\
                             &             &         & $5$-$(24,12,48)$ & \\
  Extended ternary Golay code & $[12,6,6]$ & $q=3$ & $5$-$(12,6,1)$ &  \\
  \hline
  Lifted Golay code over $\mathbb{Z}_4$  & $[24,12]$ & $\mathbb{Z}_4$ &  $5$-$(24,10,36)$ & \cite{har98,gulhar99}\\
                                        &           &                &  $5$-$(24,11,336)$ & \cite{har98}\\
                                         &           &                &  $5$-$(24,12,1584)$ & \cite{har98}\\
                               &           &            &  $5$-$(24,12,1632)$ & \cite{har98}\\
  \hline
  Extended quadric residue codes & $[24,12,9]$ & $q=3$  & $5$-$(24,9,6)$ & \cite{assma69,assma74} \\
                                 &             &        & $5$-$(24,12,576)$ & \\
                                 &             &        & $5$-$(24,15,8580)$ & \\
                                 & $[30,15,12]$ & $q=4$ & $5$-$(30,12,220)$ & \cite{assma69,assma74}\\
                                 &              &       & $5$-$(30,14,5390)$ & \\
                                 &              &       & $5$-$(30,16,123000)$ & \\
                                 & $[48,24,12]$ & $q=2$ & $5$-$(48,12,8)$ & \cite{assma69,assma74}\\
                                 &              &       & $5$-$(48,16,1365)$ & \\
                                 &              &       & $5$-$(48,20,36176)$ & \\
                                 &              &       & $5$-$(48,24,190680)$ & \\
                                 & $[48,24,15]$ & $q=3$ & $5$-$(48,12,364)$ & \cite{assma69,assma74}\\
                                 &              &       & $5$-$(48,18,50456)$ & \\
                                 &              &       & $5$-$(48,21,2957388)$ & \\
                                 &              &       & $5$-$(48,24,71307600)$ & \\
                                 &              &       & $5$-$(48,27,749999640)$ & \\
                                 & $[60,30,18]$ & $q=3$ & $5$-$(60,18,3060)$ & \cite{assma69,assma74}\\
                                 &              &       & $5$-$(60,21,449820)$ & \\
                                 &              &       & $5$-$(60,24,34337160)$ & \\
                                 &              &       & $5$-$(60,27,1271766600)$ & \\
                                 &              &       & $5$-$(60,30,24140500956)$ & \\
                                 &              &       & $5$-$(60,33,239329029060)$ & \\
\hline
Pless symmetry codes             & $[24,12,9]$  & $q=3$ & $5$-$(24,9,6)$ & \cite{pless72} \\
                                 &              &       & $5$-$(24,12,576)$ & \\
                                 &              &       & $5$-$(24,15,8580)$ & \\
                                 & $[36,18,12]$ & $q=3$ & $5$-$(36,12,45)$ & \cite{pless72}\\
                                 &              &       & $5$-$(36,15,5577)$ & \\
                                 &              &       & $5$-$(36,18,209685)$ & \\
                                 &              &       & $5$-$(36,21,2438973)$ & \\
                                 & $[48,24,15]$ & $q=3$ & $5$-$(48,12,364)$ & \cite{pless72}\\
                                 &              &       & $5$-$(48,18,50456)$ & \\
                                 &              &       & $5$-$(48,21,2957388)$ & \\
                                 &              &       & $5$-$(48,24,71307600)$ & \\
                                 &              &       & $5$-$(48,27,749999640)$ & \\
                                 & $[60,30,18]$ & $q=3$ & $5$-$(60,18,3060)$ & \cite{pless70,pless72}\\
                                 &              &       & $5$-$(60,21,449820)$ & \\
                                 &              &       & $5$-$(60,24,34337160)$ & \\
                                 &              &       & $5$-$(60,27,1271766600)$ & \\
                                 &              &       & $5$-$(60,30,24140500956)$ & \\
                                 &              &       & $5$-$(60,33,239329029060)$ & \\
  \hline
\end{tabular}
\label{table:assmatt}
\end{InTextTable}
\end{remark}

\begin{note}
The lifted Golay code over $\mathbb{Z}_4$ is defined
in~\cite{cald95} as the extended Hensel lifted quadric residue code
of length $24$. The supports of the codewords of Hamming weight $10$
in the lifted Golay code and certain extremal double circulant Type
II codes of length $24$ yield (non-isomorphic)
\mbox{$5$-$(24,10,36)$} designs. We further note that the quadratic
residue codes and the Pless symmetry codes listed in the table with
the same parameters are not equivalent as shown in~\cite{pless72} by
inspecting specific elements of the automorphism group $PSL(2,q)$.
\end{note}

\begin{remark}
The concept of the weight enumerator can be generalized to
non-linear codes (so-called \emph{distance enumerator}), which leads
to an analog of the MacWilliams relations as well as to similar
results to the Assmus-Mattson Theorem for non-linear codes
(see~\cite{will72,dels73,dels73b} and Subsection~\ref{schemes}). The
question whether there is an analogous result to the Assmus-Mattson
theorem for codes over $\mathbb{Z}_4$ proposed in~\cite{har98} was
answered in the affirmative in~\cite{tan03}. Further generalizations
of the Assmus-Mattson Theorem are known, see in
particular~\cite{cald91,kepl94,sim95,bach99,kim03,shin04,brsh08}.
\end{remark}

%%% ----------------------------------------------------------------------

\subsection{Codes and finite geometries}\label{geom}

Let $A$ be an incidence matrix of a projective plane $PG(2,n)$ of
order $n$. When we consider the subspace $C$ of
$\mathbb{F}_2^{n^2+n+1}$ spanned by the rows of $A$, we obtain for
odd $n$ only the $[n^2+n+1,n^2+n,2]$ code consisting of all
codewords of even weight. The case for even $n$ is more interesting,
in particular if $n \equiv$ $2$ (mod $4$).

\begin{theorem}\label{projpl_1}
For \mbox{$n \equiv 2$ $($\emph{mod} $4)$}, the rows of an incidence
matrix of a projective plane $PG(2,n)$ of order $n$ generate a
binary code $C$ of dimension $(n^2+n+2)/2$, and the extended code
$\overline{C}$ is self-dual.
\end{theorem}

In a projective plane $PG(2,n)$ of even order $n$, there exist sets
of $n+2$ points, no three of which are collinear, and which are
called \emph{hyperovals}\index{hyperoval} (sometimes just
\emph{ovals}\index{oval}, cf.~\cite{hirsch98}). This gives
furthermore

\begin{theorem}\label{projpl_2}
The code $C$ has minimum weight $n+1$. Moreover, the codewords of
minimum weight correspond to the lines and those of weight $n+2$ to
the hyperovals of $PG(2,n)$.
\end{theorem}

\begin{remark}\label{projpl10}
The above two theorems arose in the context of the examination of
the existence of a projective plane of order $10$
(cf.~Remark~\ref{projpl_existence}; for detailed proofs see,
e.g.,~\cite[Chapt.\,13]{cali91}). Assuming the existence of such a
plane, the obtained properties of the corresponding code lead to
very extensive computer searches. For example, in an early crucial
step, it was shown~\cite{will73} that this code could not have
codewords of weight $15$. On the various attempts to attack the
problem and the final verification of the non-existence, we refer
to~\cite{lam89,lam91,cam95} as well as~\cite[Chap.\,17]{hall86}
and~\cite[Chap.\,12]{kaos06}.
\end{remark}

\begin{note}
We note that at present the Fano plane is the only known projective
plane with order $n \equiv 2$ (mod $4)$.
\end{note}

For further accounts on codes and finite geometries, the reader is
referred, e.g., to~\cite[Chap.\,5\,\mbox{and}\,6]{assk93}
and~\cite{berl68,assk96,assk98,cam95,cali91,thas98,leo03}, as well
as~\cite{her79} from a more group-theoretical perspective
and~\cite{cam99} with an emphasis on quadratic forms over
$\mathbb{F}_2$.

%%% ----------------------------------------------------------------------

\subsection{Golay codes, Mathieu-Witt designs, and Mathieu groups}\label{golay_witt_mathieu}

We highlight some of the remarkable and natural interrelations
between the Golay codes, the Mathieu-Witt designs, and the Mathieu
groups.

\smallskip

There are various different constructions for the Golay codes
besides the description as quadratic residue codes in
Example~\ref{golay_QR}. We briefly illustrate some exemplary
constructions. For further details and more constructions, we refer
to~\cite{beju81},~\cite[Chap.\,20]{wisl77},~\cite[Chap.\,11]{cali91},
and~\cite[Chap.\,11]{cosl98}.

\begin{construction}\label{golay_constr}
\begin{itemize}

\item Starting with the zero vector in $\mathbb{F}_2^{24}$, a linear code
of length $24$ can be obtained by successively taking the
lexicographically least binary codeword which has not been used and
which has distance at least $8$ to any predecessor. At the end of
this process, we have $4096$ codewords which form the extended
binary Golay code. This construction is due to J.~H.~Conway and
N.~J.~A.~Sloane~\cite{cosl86}.

\item Let $A$ be an incidence matrix of the (unique)
\mbox{$2$-$(11,6,3)$} design. Then $G:=(I_{12},P)$ with
\[P:=\left(
  \begin{array}{cccc}
    0 & \; 1 & \;\cdots & 1 \\
    1 &  &  &  \\
    \vdots &  & A &  \\
    1 &  &  &  \\
  \end{array}
\right)\] is a $(12 \times 24)$-matrix in which each row (except the
top row) has eight $1$'s, and generates the extended binary Golay
code.

\item Let $N$ be an $(12 \times 12)$-adjacency matrix of the graph
formed by the vertices and edges of the regular icosahedron. Then
$G:=(I_{12},J_{12}-N)$ is a generator matrix for the extended binary
Golay code.

\item We recall that $\mathbb{F}_4=\{0,1,\omega, \omega^2\}$ is the field
of four elements with $\omega^2 = \omega +1$. The
\emph{hexacode}\index{hexacode} is the $[6,3,4]$ code over
$\mathbb{F}_4$ generated by the matrix $G:=(I_3,P)$ with
\[P:=\left(
  \begin{array}{ccc}
    1 & {\omega}^2 & \omega \\
    1 & \omega & {\omega}^2 \\
    1 & 1 & 1\\
  \end{array}
\right).\] The extended binary Golay code can be defined by
identifying each codeword with a binary $(4 \times 6)$-matrix $M$
(with rows $\mathbf{m}_0,\mathbf{m}_1,\mathbf{m}_2,\mathbf{m}_3$),
satisfying
\begin{enumerate}
\item[(i)] each column of $M$ has the same parity as the first row $\mathbf{m}_0$,

\item[(ii)] the sum $\mathbf{m}_1 + \omega \mathbf{m}_2 + \omega^2 \mathbf{m}_3$ lies in the hexacode.
\end{enumerate}
This description is essentially equivalent to the computational tool
\emph{MOG (Miracle Octad Generator)}\index{Miracle Octad Generator}
of R.~T.~Curtis~\cite{curt76}. The construction via the hexacode is
by Conway, see, e.g.,~\cite[Chap.\,11]{cosl98}.

\item Let $Q$ be the circulant matrix of order $5$ defined
by Eq.~(\ref{circ}). Then $G:=(I_{6},P)$, where $P$ is the matrix
$Q$ bordered on top with a row of $1$'s, is a generator matrix of
the ternary Golay code.
\end{itemize}
\end{construction}

\begin{remark}
Referring to Example~\ref{Witt_des}, we note that the automorphism
groups of the Golay codes are isomorphic to the particular Mathieu
groups, as was first pointed out in~\cite{pai56,assma66}. Moreover,
the Golay codes are related in a particularly deep and interesting
way to a larger family of sporadic finite simple groups (cf.,
e.g.,~\cite{asch94}).
\end{remark}

\begin{remark}
We have seen in Example~\ref{golay} that the supports of the
codewords of weight $8$ in the extended binary $[24,12,8]$ Golay
code form a Steiner \mbox{$5$-$(24,8,1)$} design. The uniqueness of
the large Mathieu-Witt design (up to isomorphism) can be established
easily via coding theory (cf.~Example~\ref{Witt_des}). The main part
is to show that any binary code of $4096$ codewords, including the
zero vector, of length $24$ and minimum distance $8$, is linear and
can be determined uniquely (up to equivalence). For further details,
in particular for a uniqueness proof of the small Mathieu-Witt
designs, see, e.g.,~\cite{lint93,beju81}
and~\cite[Chap.\,11]{cali91}.
\end{remark}

%%% ----------------------------------------------------------------------

\subsection{Golay codes, Leech lattice, kissing numbers, and sphere packings}\label{golay_lattice}

Sphere packings closely connect mathematics and information
theory\index{information theory} via the sampling
theorem\index{Shannon's Sampling Theorem} as observed by
C.~E.~Shannon~\cite{shan48} in his classical article of 1948.
Rephrased in a more geometric language, this can be expressed as
follows:

\begin{quote}
``Nearly equal signals are represented by neighboring points, so to
keep the signals distinct, Shannon represents them by
$n$-dimensional `billiard balls', and is therefore led to ask: what
is the best way to pack `billiard balls' in $n$
dimensions?''\cite{slo98}
\end{quote}

One of the most remarkable lattices, the \emph{Leech
lattice}\index{Leech lattice} in $\mathbb{R}^{24}$, plays a crucial
role in classical sphere packings. We recall that a
\emph{lattice}\index{lattice} in $\mathbb{R}^n$ is a discrete
subgroup of $\mathbb{R}^n$ of rank $n$. The extended binary Golay
code led to the discovery by John Leech~\cite{leech64} of the
\mbox{$24$-dimensional} Euclidean lattice named after him. There are
various constructions besides the usual ones from the binary and
ternary Golay codes in the meantime, see,
e.g.,~\cite{cosl82},~\cite[Chap.\,24]{cosl98}. We outline some of
the fundamental connections between sphere packings and the Leech
lattice.

\smallskip

The \emph{Kissing Number Problem}\index{Kissing Number Problem}
deals with the maximal number $\tau_n$ of equal size non-overlapping
spheres in the $n$-dimensional Euclidean space $\mathbb{R}^n$ that
can touch a given sphere of the same size. Only a few of these
numbers are actually known. For dimensions $n=1,2,3$, the classical
solutions are: $\tau_1=2$, $\tau_2=6$, $\tau_3=12$. The number
$\tau_3$ was the subject of a famous controversy between Isaac
Newton and David Gregory in 1694, and was finally verified only in
1953 by K.~Sch{\"u}tte and B.~L.~van der Waerden~\cite{sch53}. Using
an approach initiated by P.~Delsarte~\cite{dels72,dels73} in the
early 1970's which gives linear programming upper bounds for binary
error-correcting codes and for spherical codes~\cite{dels77}
(cf.~Subsection~\ref{schemes}), A.~M.~Odlyzko and
N.~J.~A.~Sloane~\cite{odslo79}, and independently
V.~I.~Levenshtein~\cite{lev79}, proved that $\tau_8 = 240$ and
$\tau_{24} = 196560$. These exact solutions are the number of
non-zero vectors of minimal length in the root lattice $E_8$ and in
the Leech lattice, respectively. By extending and improving
Delsarte's method, O.~R.~Musin~\cite{musin} verified in 2003 that
$\tau_4=24$, which is the number of non-zero vectors of minimal
length in the root lattice $D_4$.

\smallskip

The \emph{Sphere Packing Problem}\index{Sphere Packing Problem} asks
for the maximal density of a packing of equal size non-overlapping
spheres in the \mbox{$n$-dimensional} Euclidean space
$\mathbb{R}^n$. A sphere packing is called a \emph{lattice
packing}\index{lattice packing} if the centers of the spheres form a
lattice in $\mathbb{R}^n$. The Leech lattice is the unique densest
lattice packing (up to scaling and isometries) in $\mathbb{R}^{24}$,
as was shown by H.~Cohn and A.~Kumar~\cite{cokum04,cokum04b}
recently in 2004, again by a modification of Delsarte's method.
Moreover, they showed that the density of any sphere packing in
$\mathbb{R}^{24}$ cannot exceed the one given by the Leech lattice
by a factor of more than $1 + 1.65 \cdot 10^{-30}$ (via a computer
calculation). The proof is based on the work~\cite{coelk03} by Cohn
and N.~D.~Elkies in 2003 in which linear programming bounds for the
Sphere Packing Problem are introduced and new upper bounds on the
density of sphere packings in $\mathbb{R}^{n}$ with dimension $n
\leq 36$ are proven.

\smallskip

For further details on the Kissing Number Problem and the Sphere
Packing Problem,
see~\cite[Chap.\,1]{cosl98},~\cite{musin05,slo98},~\cite{thomp83},
as well as the survey articles~\cite{elk00,zieg04,bez06}. For an
on-line database on lattices, see~\cite{nesl}.

%%% ----------------------------------------------------------------------

\subsection{Codes and association schemes}\label{schemes}

Any finite nonempty subset of the unit sphere $S^{n-1}$ in the
$n$-dimensional Euclidean space $\mathbb{R}^n$ is called a
\emph{spherical code}\index{spherical code}. These codes have many
practical applications, e.g., in the design of signals for data
transmission and storage. As a special class of spherical codes,
\emph{spherical designs}\index{spherical design} were introduced by
P.~Delsarte, {J.-M.}~Goethals and J.~Seidel~\cite{dels77} in 1977 as
analogs on $S^{n-1}$ of the classical combinatorial designs. For
example, in $S^2$ the tetrahedron is a spherical \mbox{$2$-design};
the octahedron and the cube are spherical \mbox{$3$-designs}, and
the icosahedron and the dodecahedron are spherical
\mbox{$5$-designs}. In order to obtain the linear programming upper
bound mentioned in the previous subsection, Krawtchouk polynomials
were involved in the case of binary error-correcting codes and
Gegenbauer polynomials in the case of spherical codes.

\smallskip

However, Delsarte's approach was indeed much more general and
far-reaching. He developed for association schemes, which have their
origin in the statistical theory of design of experiments, a theory
to unify many of the objects we have been addressing in this
chapter. We give a formal definition of \emph{association
schemes}\index{association scheme} in the sense of
Delsarte~\cite{dels73} as well as introduce the
\emph{Hamming}\index{Hamming scheme} and the \emph{Johnson
schemes}\index{Johnson scheme} as important examples of the two
fundamental classes of \emph{\mbox{$P$-polynomial}} and
\emph{\mbox{$Q$-polynomial} association schemes}.

\begin{definition}\label{ass_schemes}
A \emph{$d$-class association scheme} is a finite point set $X$
together with $d+1$ relations $R_i$ on $X$, satisfying
\begin{enumerate}
\item[(i)] $\{R_0,R_1,\ldots,R_d\}$ is a partition of $X \times X$,

\item[(ii)] $R_0=\{(x,x)\mid x \in X\}$,

\item[(iii)] for each $i$ with $0 \leq i \leq d$, there exists a $j$ with $0 \leq j \leq d$ such that $(x,y) \in R_i$ implies $(y,x) \in R_j$,

\item[(iv)] for any  $(x,y) \in R_k$, the number $p_{ij}^k$ of points $z \in X$ with $(x,z) \in R_i$ and $(z,y) \in R_j$
depends only on $i,j$ and $k$,

\item[(v)] $p_{ij}^k=p_{ji}^k$ for all $i$, $j$ and $k$.
\end{enumerate}
\end{definition}
The numbers $p_{ij}^k$ are called the \emph{intersection numbers} of
the association scheme. Two points $x,y \in X$ are called
\emph{\mbox{$i$-th} associates} if $\{x,y\} \in R_i$.

\begin{example}
The \emph{Hamming scheme} $H(n,q)$ has as point set $X$ the set
$\mathbb{F}^n$ of all \mbox{$n$-tuples} from a \mbox{$q$-symbol}
alphabet; two $n$-tuples are \mbox{$i$-th} associates if their
Hamming distance is $i$. The \emph{Johnson scheme} $J(v,k)$, with $k
\leq \frac{1}{2}v$, has as point set $X$ the set of all
\mbox{$k$-element} subsets of a set of size $v$; two
\mbox{$k$-element} subset $S_1,S_2$ are \mbox{$i$-th} associates if
$\left|S_1 \cap S_2 \right|=k-i$.
\end{example}

Delsarte introduced the Hamming and Johnson schemes as settings for
the classical concept of error-correcting codes and combinatorial
\mbox{$t$-designs}, respectively. In this manner, certain results
become formally dual, like the Sphere Packing Bound
(Theorem~\ref{Hamming-bound}) and Fisher's Inequality
(Theorem~\ref{FisherIn}).

\smallskip

For a more extended treatment of association schemes, the reader is
referred, e.g.,
to~\cite{bos52,bos59,bait84,bcn89,brou95},~\cite[Chap.\,17]{cali91},~\cite[Chap.\,21]{wisl77},
and in particular to~\cite{cam98,dels98} with an emphasis on the
close connection between coding theory and associations schemes. For
a survey on spherical designs, see~\cite[Chap.\,VI.54]{crc06}.

%%% ----------------------------------------------------------------------

\section{Directions for further research}\label{directions}

We present in this section a collection of significant open problems
and challenges concerning future research.

\begin{problem}\em{(cf.~\cite{stein53}).}
Does every Steiner triple system on $n$ points extend to a Steiner
quadruple system on $n + 1$ points?
\end{problem}

\begin{problem}\em
Does there exist any non-trivial Steiner \mbox{$6$-design}?
\end{problem}

\begin{problem}\em{(cf.~\cite[p.\,180]{wisl77}).}
Find all non-linear single-error-correcting perfect codes over
$\mathbb{F}_q$.
\end{problem}

\begin{problem}\em{(cf.~\cite[p.\,106]{hupl_handb98}).}
Characterize codes where all codewords of the same weight (or of
minimum weight) form a non-trivial design.
\end{problem}

\begin{problem}\em{(cf.~\cite[p.\,116]{hupl_handb98}).}
Find a proof of the non-existence of a projective plane of order
$10$ without the help of a computer or with an easily reproducible
computer program.
\end{problem}

\begin{problem}\em
Does there exist any finite projective plane of order $12$, or of
any other order that is neither a prime power nor covered by the
Bruck-Ryser Theorem (cf.~Remark~\ref{projpl_existence})?
\end{problem}

\begin{problem}\em
Does the root lattice $D_4$ give the unique kissing number
configuration in $\mathbb{R}^{4}$?
\end{problem}

\begin{problem}\em
Solve the  Kissing Number Problem in $n$ dimensions for any $n
> 4$ apart from $n=8$ and $24$. For presently known lower and upper bounds, we refer
to~\cite{kissnb} and~\cite{val07}, respectively. Also any improvements of these bounds
would be desirable.
\end{problem}

\begin{problem}\em{(cf.~\cite[Conj.\,8.1]{coelk03}).}
Verify the conjecture that the Leech lattice is the unique densest
sphere packing in $\mathbb{R}^{24}$.
\end{problem}

%%% ----------------------------------------------------------------------

\section{Conclusions}\label{conclusions}

Over the last sixty years a substantial amount of research has been
inspired by the various interactions of coding theory and algebraic
combinatorics. The fruitful interplay often reveals the high degree
of regularity of both the codes and the combinatorial structures.
This has lead to a vivid area of research connecting closely
mathematics with information and coding theory. The emerging methods
can be applied sometimes surprisingly effectively, e.g., in view of
the recent advances on kissing numbers and sphere packings.

\smallskip

A further development of this beautiful interplay as well as its
application to concrete problems would be desirable, certainly also
in view of the various still open and long-standing problems.

%%% ----------------------------------------------------------------------

\section{Terminologies/Keywords}

Error-correcting codes, combinatorial designs, perfect codes and
related concepts, Assmus-Mattson Theorem and analogues, projective
geometries, non-existence of a projective plane of order $10$, Golay
codes, Leech lattice, kissing numbers, sphere packings, spherical
codes, association schemes.

%%% ----------------------------------------------------------------------

\section{Exercises}\label{exercises}

\begin{enumerate}

\item[(1)] Verify (numerically) that the Steiner quadruple system $SQS(8)$ of order $8$ (cf.~Example~\ref{AG_cube}) has $14$ blocks,
and that the Mathieu-Witt design \mbox{$5$-$(24,8,1)$}
(cf.~Example~\ref{Witt_des}) has $759$ blocks.

\smallskip

\item[(2)] What are the parameters of the \mbox{$2$-design} consisting of the points and hyperplanes (i.e. the $(d-2)$-dimensional
projective subspaces) of the projective space $PG(d - 1, q)$?

\smallskip

\item[(3)] Does there exist a self-dual $[8,4]$ code over the finite field $\mathbb{F}_2$?

\smallskip

\item[(4)] Show that the ternary $[11,6,5]$ Golay code has $132$
codewords of weight $5$.

\smallskip

\item[(5)] Compute the weight distribution of the binary $[23,12,7]$ Golay code.

\smallskip

\item[(6)] Show that any binary code of $4096$ codewords, including the zero
vector, of length $24$ and minimum distance $8$ is linear.

\smallskip

\item[(7)] Give a proof for the Sphere Packing Bound (cf.~Theorem~\ref{Hamming-bound}).

\smallskip

\item[(8)] Give a proof for Fisher's Inequality (cf.~Theorem~\ref{FisherIn}).
\smallskip

\item[(9)] Show that a binary code generated by the rows of an incidence matrix of any projective plane $PG(2,n)$
of even order $n$ has dimension at most $(n^2 + n + 2)/2$
(cf.~Theorem~\ref{projpl_1}).

\smallskip

\item[(10)] (Todd's Lemma). In the Mathieu-Witt design
\mbox{$5$-$(24,8,1)$}, if $B_1$ and $B_2$ are blocks
(\emph{octads}\index{octad}) meeting in four points, then $B_1 +
B_2$ is also a block.

\end{enumerate}

\bigskip

\textbf{Solutions:}

\smallskip

\begin{enumerate}

\item[ad (1):] By Lemma~\ref{Comb_t=5}~(b), we have to calculate $b=\frac{8 \cdot 7 \cdot
6}{4 \cdot 3 \cdot 2}=14$ in the case of the Steiner quadruple
system $SQS(8)$, and $b=\frac{24 \cdot 23 \cdot 22 \cdot 21 \cdot
20}{8 \cdot 7 \cdot 6 \cdot 5 \cdot 4}=759$ in the case of the
Mathieu-Witt design \mbox{$5$-$(24,8,1)$}.

\smallskip

\item[ad (2):] Starting from Example~\ref{PG}, we obtain via counting arguments (or by using
the transitivity properties of the general linear group) that the
points and hyperplanes of the projective space $PG(d - 1, q)$ form a
\mbox{$2$-$(\frac{q^{d}-1}{q-1},\frac{q^{d-1}-1}{q-1},\frac{q^{d-2}-1}{q-1})$
design}.

\smallskip

\item[ad (3):] Yes, the extended binary $[8,4,4]$ Hamming code is self-dual
(cf.~Example~\ref{H8}).

\smallskip

\item[ad (4):] Since the ternary $[11,6,5]$ Golay code is perfect (cf.~Example~\ref{golay}), every word of
weight $3$ in $\mathbb{F}_3^{11}$ has distance $2$ to a codeword of
weight $5$. Thus $A_5= 2^3 \cdot {11 \choose 3}/{5 \choose 2}=132$.

\smallskip

\item[ad (5):] The binary $[23,12,7]$ Golay code contains the zero vector and is perfect. This
determines the weight distribution as follows $A_0 =A_{23}=1$,
$A_7=A_{16}=253$, $A_8=A_{15}=506$, $A_{11}=A_{12}=1288$.

\smallskip

\item[ad (6):] Let $C$ denote a binary code of $4096$ codewords, including the zero
vector, of length $24$ and minimum distance $8$. Deleting any
coordinate leads to a code which has the same weight distribution as
the code given in Exercise~(5). Hence, the code $C$ only has
codewords of weight $0,8,12,16$ and $24$. This is still true if the
code $C$ is translated by any codeword (i.e. $C + \mathbf{x}$ for
any $\mathbf{x} \in C$). Thus, the distances between pairs of
codewords are also $0, 8, 12, 16$ and $24$. Therefore, the standard
inner product $\langle \textbf{x},\textbf{y}\rangle$ vanishes for
any two codewords $\textbf{x},\textbf{y} \in C$, and hence $C$ is
self-orthogonal. For cardinality reasons, we conclude that $C$ is
self-dual and hence in particular linear.

\smallskip

\item[ad (7):] The sum $\sum_{i=0}^{e} {n \choose i} (q-1)^i$ counts
the number of words in a sphere of radius $e$. As the spheres of
radius $e$ about distinct codewords are disjoint, we obtain $\left|
C \right| \cdot \sum_{i=0}^{e} {n \choose i} (q-1)^i$ words.
Clearly, this number cannot exceed the total number $q^n$ of words,
and the claim follows.

\smallskip

\item[ad (8):] As a non-trivial \mbox{$t$-design} with $t \geq 2$ is also a
non-trivial \mbox{$2$-design} by Lemma~\ref{s-design}, it is
sufficient to prove the assertion for an arbitrary non-trivial
\mbox{$2$-$(v,k,\lambda)$} design $\mathcal{D}$. Let $A$ be an
incidence matrix of $\mathcal{D}$ as defined in
Subsection~\ref{intro_des}. Clearly, the $(i,k)$-th entry
\[(AA^t)_{ik}= \sum_{j=1}^{b} (A)_{ij}(A^t)_{jk}=\sum_{j=1}^{b} a_{ij}a_{kj}\]
of the $(v \times v)$-matrix $AA^t$ is the total number of blocks
containing both $x_i$ and $x_k$, and is thus equal to $r$ if $i=k$,
and to $\lambda$ if $i \not=k$. Hence
\[AA^t=(r- \lambda) I + \lambda J,\]
where $I$ denotes the $(v \times v)$-unit matrix and $J$ the $(v
\times v)$-matrix with all entries equal to $1$. Using elementary
row and column operations, it follows easily that
\[\det(AA^t)=r k  (r - \lambda)^{v-1}.\] Thus $AA^t$ is
non-singular (i.e. its determinant is non-zero) as $r=\lambda$ would
imply $v=k$ by Lemma~\ref{s-design}, yielding that the design is
trivial. Therefore, the matrix $AA^t$ has rank$(A)=v$. But, if
$b<v$, then rank$(A) \leq b <v$, and thus rank$(AA^t)<v$, a
contradiction. It follows that $b \geq v$, proving the claim.

\smallskip

\item[ad (9):] Let $C$ denote a binary code generated by the rows of an incidence matrix of $PG(2,n)$.
By assumption $n$ is even, and hence the extended code
$\overline{C}$ must be self-orthogonal. Therefore, the dimension of
$C$ is at most $n^2+n+2/2$.

\smallskip

\item[ad (10):] For given blocks $B_1=\{01,02,03,04,05,06,07,08\}$ and
$B_2=\{01,02,03,04,\linebreak 09,10,11,12\}$ in the Mathieu-Witt
design \mbox{$5$-$(24,8,1)$}, let us assume that $B_1+B_2$ is not a
block. The block $B_3$ which contains $\{05,06,07,08,09\}$ must
contain just one more point of $B_2$, say
$B_3=\{05,06,07,08,09,10,13,14\}$. Similarly,
$B_4=\{05,06,07,08,11,12,15,16\}$ is the block containing
$\{05,06,07,08,11\}$. But hence, it is impossible to find a block
which contains $\{05,06,07,09,11\}$ and intersects with $B_i$, $1
\leq i \leq 4$, in $0$, $2$ or $4$ points. Therefore, we obtain a
contradiction as there must be a block containing any five points by
Definition~\ref{StDes}.
\end{enumerate}

\bibliographystyle{ws-rv-van}
\bibliography{ws-rv-Huber_CodingAlgComb}

%\printindex[aindx]                 % to print author index
\printindex                         % to print subject index
\end{document}